\documentclass[12pt,amssymb]{article}
\usepackage{amsfonts,latexsym}

\newcommand{\proj}{\mathbb P}
\newcommand{\reel}{\mathbb R}
\newcommand{\complex}{\mathbb C}

\newtheorem{theorem}{Th\'{e}or\`{e}me}
\newtheorem{lem}{Lemme}

\begin{document}
\title{Enveloppes inf\'{e}rieures de fonctions admissibles
sur l'espace projectif complexe. Cas sym\'etrique.}
\author {\em A. Ben Abdesselem}
\date{\em Universit\'e Pierre et Marie Curie, U.F.R. 920, BC 172, 4, place Jussieu, 75005, Paris, France}

\maketitle {\em RESUME. On montre que les fonctions admissibles
pour la m\'etrique de Fubini-Study sur l'espace projectif complexe
$\proj_{m}\complex$ de dimension complexe $m$, invariantes par un
groupe d'automorphismes convenablement choisi, sont minor\'ees par
une fonction tendant vers moins l'infini sur le bord des cartes
usuelles de $\proj_{m}\complex$. Une minoration similaire a lieu
sur certaines vari\'et\'es projectives. Cette minoration donne une
constante optimale dans une in\'egalit\'e de type H\H{o}rmander
sur ces vari\'et\'es, ce qui permet d'y \'etablir l'existence de
m\'etriques d'Einstein-K\"{a}hler.}

\maketitle {\em ABSTRACT. We prove that admissible functions for
Fubini-Study metric on the complex projective space
$\proj_{m}\complex$ of complex dimension $m$, invariant by a
convenient automorphisms group, are  lower bounded by a function
going to minus infinity on the boundary of usual charts of
$\proj_{m}\complex$. A similar lower bound holds on some
projective manifolds. This gives an optimal constant in a
H\H{o}rmander type inequality on these manifolds, which allows us
to establish the existence of Einstein-K\"{a}hler metrics on them.

}

{\section{Introduction.}\label{i}}

\mbox{\,\,\,\,\,\,\,}

On consid\`{e}re l'espace projectif complexe $\proj_{m}\complex$
de dimension complexe $m=kn-1$.  Si $[z_{0},z_{1},..,z_{m}]$
d\'esignent les coordonn\'{e}es homog\`{e}nes de
$\proj_{m}\complex$, on d\'efinit, pour $h\in \{0,..,k-1\}$, le
$n$-uplet $Z_{h}\in \complex^{n}$ par $Z_{h}=
(z_{hn},..,z_{(h+1)n-1})$. Ceci nous permet d'\'ecrire
$[z_{0},z_{1},..,z_{m}]=[Z_{0},..,Z_{k-1}]\in \proj_{m}\complex$.
D'autre part, on munit $\proj_{m}\complex$ de la m\'{e}trique $g$
ayant pour composantes, dans la carte $\{z_{0}\neq 0\}$,
$$g_{\lambda\overline{\mu}}=a_{m}\frac{\partial^{2}}{\partial
z_{\lambda}\partial\overline{z}_{\mu}}(\ln (1+\mid z_{1} \mid
^{2}+..+\mid z_{m} \mid ^{2})),$$ o\`u $a_{m}>0$. Elle est
\'{e}gale \`{a} $a_{m}$ fois la m\'{e}trique de Fubini-Study et
appartient \`{a} la premi\`{e}re classe de Chern
$C_{1}(\proj_{m}\complex )$ lorsque $a_{m}=(m+1)$. On dit que
$\varphi\in C^{\infty}(\proj_{m}\complex)$ est $g$-admissible
lorsque
$$g_{\lambda\overline{\mu}}+\frac{\partial^{2}\varphi}{\partial
z_{\lambda}\partial\overline{z}_{\mu}}
>0.$$
On consid\`{e}re aussi le groupe d'automorphismes $G_{n,k}$ sur
$\proj_{m}\complex$ engendr\'{e} par les automorphismes
$$\sigma_{i,j}:[Z_{0},..,Z_{i},..,Z_{j},..,Z_{k-1}]
\longrightarrow[Z_{0},..,Z_{j},..,Z_{i},..,Z_{k-1}],$$
$$\tau_{k,\theta}:[z_{0},...,z_{m}]
\longrightarrow[z_{0},..,z_{k}e^{i\theta},..,z_{m}],$$ et enfin,
pour $z_{p}$ et $z_{q}$ appartenant \`a un m\^eme $n$-uplet
$Z_{h}$, $$\gamma_{p,q}:[z_{0},..,z_{p},..,z_{q},..,z_{m}]
\longrightarrow[z_{0},..,z_{q},..,z_{p},..,z_{m}].$$

On d\'{e}finit sur $\complex^{m+1}\backslash\bigcup_{p}\{z_{p}=0\}$ la
fonction $\psi$ par $$\psi =\ln \frac{(\mid z_{0}\mid...\mid
z_{m}\mid )^{2a_{m}/(m+1)}}{(\mid z_{0}\mid^{2}+...+\mid
z_{m}\mid^{2})^{a_{m}}}.$$ Elle est homog\`{e}ne de degr\'{e}
z\'{e}ro sur $\complex^{m+1}$, elle d\'{e}finit donc une fonction
sur $\proj_{m}\complex$. $\psi$ atteint son maximum (\'egal \`a
$-a_{m}\ln (m+1)$) en les points
$$[1,e^{i\theta_{1}},..,e^{i\theta_{m}}] \in\proj_{m}\complex$$
qui d\'{e}finissent un produit de tores. $\psi$ tends vers
moins l'infini lorsque l'une des coordonn\'{e}es homog\`{e}nes
$z_{0}$,.., $z_{m}$ tend vers z\'{e}ro ou vers l'infini, ce qui
correspond aux fronti\`{e}res des cartes denses definies par
$\{z_{h}=0\}$.

Etablissons \`a pr\'esent les principaux r\'{e}sultats de cet
article :
\begin{theorem}\label{th1}.
Soit $\varphi\in C^{\infty}(\proj_{m}\complex)$ une fonction
$g$-admissible et $G_{n,k}$-invariante, v\'{e}rifiant $\sup\varphi
= 0$ sur $\proj_{m}\complex$. On a alors: $\varphi \geq \psi$.
\end{theorem}
On a alors le corollaire suivant.
\begin{theorem}\label{th2}.
Pour tout $\alpha < 1$, on a l'in\'{e}galit\'{e} de type
H\H{o}rmander (voir \cite{H} th. 4.4.5):
$$\int_{\proj_{m}\complex} \exp\{-\alpha\varphi \} dv \leq Cst,$$
$\forall \varphi\in C^{\infty}(\proj_{m}\complex)$, $g$-admissible
avec $a_{m}=(m+1)$, $G_{n,k}$-invariante, et v\'{e}rifiant $\sup
\varphi =0$ sur $\proj_{m}\complex$. $dv$ est l'\'{e}l\`{e}ment de
volume sur $\proj_{m}$ relatif \`{a} la m\'{e}trique $g$.
\end{theorem}


L'in\'{e}galit\'{e} du th\'{e}or\`{e}me \ref{th2} est optimale car
elle correspond \`a la valeur de l'invariant de Tian (cf.
\cite{T}) relatif \`{a} cette classe de fonctions (voir par
exemple \cite{A2} ou \cite{BC1}). Cet article donne alors un moyen
simple de calculer cet invariant en r\'{e}duisant la classe de
fonctions \`{a} consid\'{e}rer \`{a} la fonction $\psi$. La
m\'ethode employ\'ee s'applique \`a des vari\'et\'es plus
compliqu\'ees que $\proj_{m}\complex$. Pour cel\`a, il suffit de
mettre en \'evidence la bonne fonction extr\'emale. En effet,
Consid\'erons par exemple la sous-vari\'et\'e $M$ du produit
$\proj_{m}\complex\times\proj_{n-1}\complex\times
...\times\proj_{n-1}\complex$ de $\proj_{m}\complex$ par $k$
exemplaires de $\proj_{n-1}\complex$ o\`u $m=kn-1$, et
constitu\'ee des points
$$([Z_{0},..,Z_{k-1}],[\zeta_{0}],..,[\zeta_{k-1}]) \in
\proj_{m}\complex\times\proj_{n-1}\complex\times
...\times\proj_{n-1}\complex ,$$ o\`u $Z_{h}=
(z_{hn},..,z_{(h+1)n-1})\in [\zeta_{h}]$. $M$ est de dimension
complexe \'egale \`a $m$. Consid\`erons la m\'etrique
$g^{M}=\tilde{g}_{|M}$ restriction \`a $M$ de la m\'etrique
$\tilde{g}$ sur $\proj_{m}\complex\times\proj_{n-1}\complex\times
...\times\proj_{n-1}\complex$, \'egale \`a
$a_{m_{0}}g^{(m_{0})}+...+a_{m_{k}}g^{(m_{k})}$ o\`u les
$g^{(m_{i})}$ sont les m\'etriques de Fubini-Study sur
$\proj_{m}\complex$, $\proj_{n-1}\complex$,..,
$\proj_{n-1}\complex$. Il est montr\'e dans \cite{BC0} que
$\tilde{g}_{|M}$ est dans la premi\`ere classe de Chern $C_{1}(M)$
pour $a_{m_{0}}=k$ et $a_{m_{1}}=...=a_{m_{k}}=n-1$. Son
expression dans une carte o\`u les $Z_{i}$ ont au moins une
composante non nulle est donn\'ee, dans ce cas par :
\begin{eqnarray}
g^{M}_{\lambda\overline{\mu}}=\frac{\partial^{2}}{\partial
z_{\lambda}\partial \overline{z}_{\mu}}\ln [(1+\mid
z_{1}\mid^{2}+...+\mid z_{kn-1}\mid^{2})^{k}\times (1+\mid
z_{1}\mid^{2}+...+\mid z_{n-1}\mid^{2})^{n-1}\times\nonumber\\
(\mid z_{n}\mid^{2}+...+\mid z_{2n-1}\mid^{2})^{n-1}\times
...\times(\mid z_{(k-1)n}\mid^{2}+...+\mid
z_{kn-1}\mid^{2})^{n-1}],\nonumber
\end{eqnarray}
(ici, les coordonn\'ees homog\`enes sont divis\'ees par $z_{0}$).
D'autre part, les groupes d'automorphismes $G_{n,k}$ de
$\proj_{m}\mathbb{C}=\mathbb{P}_{kn-1}\mathbb{C}$ et $G_{1,n}$ de
$\mathbb{P}_{n-1}\mathbb{C}$ d\'efinis plus haut, engendrent un
groupe d'automorphismes sur le produit
$\mathbb{P}_{m}\mathbb{C}\times\mathbb{P}_{n-1}\mathbb{C}\times
...\times\mathbb{P}_{n-1}\mathbb{C}$. $M$ \'etant \`a l'\'evidence
stable par ce groupe, ce dernier induit sur $M$ un groupe
d'automorphismes not\'e $G^{M}$.

Enfin, on d\'{e}finit sur $\complex^{m+1}\times(\complex^{n})^{k}
\backslash\bigcup_{i,p}{\{z_{p}^{(i)}=0\}}$, o\`{u} $m=kn-1$, la
fonction $\hat{\psi}$ par
\begin{eqnarray}
\hat{\psi} =\ln [\frac{(\mid z_{0}^{(0)}\mid...\mid
z_{m}^{(0)}\mid )^{2k/(m+1)}}{(\mid z_{0}^{(0)}\mid^{2}+...+\mid
z_{m}^{(0)}\mid^{2})^{k}}\times \frac{(\mid z_{0}^{(1)}\mid...\mid
z_{n-1}^{(1)}\mid )^{2(n-1)/n}}{(\mid z_{0}^{(1)}\mid^{2}+...+\mid
z_{n-1}^{(1)}\mid^{2})^{(n-1)}}\times...\nonumber\\ \times
\frac{(\mid z_{0}^{(k)}\mid...\mid z_{n-1}^{(k)}\mid
)^{2(n-1)/n}}{(\mid z_{0}^{(k)}\mid^{2}+...+\mid
z_{n-1}^{(k)}\mid^{2})^{(n-1)}}].\nonumber
\end{eqnarray}
$(z_{p}^{(0)})_{p}$ \'etant le syst\`eme de coordonn\'ees usuel
sur $\complex^{m+1}$, $(z_{p}^{(1)})_{p}$ celui sur le premier
$\complex^{n}$ du produit,..., et $(z_{p}^{(k)})_{p}$ celui sur le
$k$-i\`{e}me $\complex^{n}$ (on rappelle que $m=kn-1$).

Elle est homog\`{e}ne de degr\'{e} z\'{e}ro sur chacun des
$\complex^{p}$, elle d\'{e}finit donc une fonction, toujours
not\'ee $\hat{\psi}$ sur
$\proj_{m}\complex\times\proj_{n-1}\complex\times...\times\proj_{n-1}\complex$
\'egale \`a moins l'infini d\`es que l'une des composantes
homog\`enes s'annule. Sa restriction $\psi_{M}$ \`a $M$ est
$G^{M}$-invariante et d\'efinie sur l'ouvert dense de $M$
constitu\'e par les points
\begin{eqnarray}
&&([\zeta_{0}(z_{0},..,z_{n-1}),..,\zeta_{k-1}(z_{k(n-1)},..,z_{m=kn-1})],\nonumber\\
&&[z_{0},..,z_{n-1}],..,[z_{k(n-1)},..,z_{m=kn-1}])\in M\nonumber
\end{eqnarray}
o\`{u} $[\zeta_{0},..,\zeta_{k-1}]\in \mathbb{P}_{k}\mathbb{C}$
avec tous les $z_{i}$ et $\zeta_{i}$ sont non nuls. Dans le
compl\'ementaire de ces points, elle vaut moins l'infini. Sachant
que $\frac{k}{m+1}+\frac{n-1}{n}= \frac{k}{kn}+\frac{n-1}{n}=1$,
son expression donn\'ee par :
\begin{eqnarray}
\psi_{M}=\ln \{[\mid\zeta_{0}\mid^{2}(\mid z_{0}\mid^{2}+..+\mid
z_{n-1}\mid^{2})+..+\mid\zeta_{k-1}\mid^{2}(\mid
z_{(k-1)n}\mid^{2}+..+\mid z_{kn-1}\mid^{2})]^{-k}\nonumber\\
\times\frac {\mid z_{0}\mid^{2}...\mid z_{m}\mid^{2}\mid
\zeta_{0}\mid^{2}...\mid \zeta_{k-1}\mid^{2}}{(\mid
z_{0}\mid^{2}+...+\mid z_{n-1}\mid^{2})^{(n-1)}...(\mid
z_{(k-1)n}\mid^{2}+..+\mid z_{m}\mid^{2})^{(n-1)}}\}.\nonumber
\end{eqnarray}
Et en posant $z'_{i}=\zeta_{h}z_{i}$, pour
$i\in\{hn,..,(h+1)n-1\}$
\begin{eqnarray}\label{psi}
\psi_{M}=\ln \{(\mid z'_{0}\mid^{2}+..+\mid
z'_{n-1}\mid^{2}+..+\mid z'_{(k-1)n}\mid^{2}+..+\mid
z'_{kn-1}\mid^{2})^{-k}\nonumber\\ \times\frac {\mid
z'_{0}\mid^{2}...\mid z'_{m}\mid^{2}}{(\mid
z'_{0}\mid^{2}+...+\mid z'_{n-1}\mid^{2})^{(n-1)}...(\mid
z'_{(k-1)n}\mid^{2}+..+\mid z'_{m}\mid^{2})^{(n-1)}}\}.
\end{eqnarray}
On a le r\'esultat suivant :
\begin{theorem}\label{th3}.
Soit $\varphi\in C^{\infty}(M)$ une fonction $g^{M}$-admissible et
$G^{M}$-invariante telle que $\sup\varphi =0$ sur $M$. On a alors
$\varphi\geq \psi_{M}$, et $M$ admet une m\'{e}trique
d'Einstein-K\"{a}hler dans la classe de K\"{a}hler de $g^{M}$.
\end{theorem}

La m\'ethode employ\'ee dans ce papier pr\'esente l'avantage de
rendre calculable l'invariant de Tian et par cons\'equent de
donner un crit\`ere d'existence simple de m\'{e}triques
d'Einstein-K\H{a}hler. En effet, le calcul pratique de l'invariant
de Tian peut s'av\'{e}rer compliqu\'{e} comme le prouvent les
exemples trait\'{e}s dans \cite{BC1}, \cite{BC2} et \cite{X}.
D'autre part, il est naturel de chercher le lien entre le moyen
utilis\'e dans ce papier pour estimer la convergence de
$\int_{M}e^{-\alpha\varphi}dv$ pour une famille de fonctions
admissibles \`a l'aide de la fonction $\psi$ et celui utilis\'e
\`a cette m\^eme fin par A.M. Nadel dans \cite{N}.
L'int\'{e}r\^{e}t g\'{e}om\'{e}trique d'une telle estimation est
de prouver l'existence de m\'{e}triques d'Einstein-K\H{a}hler sur
les vari\'{e}t\'{e}s que nous consid\'erons. Dans le m\^eme
domaine on citera les papiers pr\'ecurseurs de T. Aubin \cite{A1}
et S.T. Yau \cite{Y} et Y.T. Siu \cite{S}. Il est d'autre part
recommand\'e de lire \cite{A3} ou le s\'eminaire Bourbaki de J.P.
Bourguignon \cite{JPB} qui font le point sur la question.


{\section{Preuve des r\'{e}sultats \'{e}nonc\'{e}s.}\label{p}}

L'id\'ee de la d\'emonstration du th\'{e}or\`{e}me \ref{th1} est
celle utilis\'ee dans \cite{B} sur $\proj_{1}\complex$. Sa mise en
oeuvre n\'ecessite les lemmes \ref{lem1}, \ref{lem2}, \ref{lem3}
et \ref{lem4} . Le th\'{e}or\`{e}me \ref{th2} est un corollaire du
premier th\'eor\`eme. D'autre part, on utilisera l'invariance des
fonctions $\varphi ([z_{0},...,z_{m}])$ par le groupe $G_{n,k}$
pour les consid\'erer comme des fonctions $\varphi
([1,x_{1},...,x_{m}])$, des variables r\'eelles $0\leq x_{i} =\mid
z_{i} \mid \leq 1$. La m\'ethode employ\'ee dans la preuve du
th\'eor\`eme \ref{th1} s'applique pour th\'eor\`eme \ref{th3} en
changeant $\psi$ par $\psi_{M}$.


\begin{lem}\label{lem1}
Etant donn\'{e}e une fonction $\varphi\in
C^{\infty}(\proj_{m}\complex)$, $g$-admissible sur
$\proj_{m}\complex$, avec $m=kn-1$, $G_{n,k}$-invariante, on a
pour $0<x_{i}\leq 1$:
\begin{eqnarray}\label{eq1}
(\varphi -\psi )([1,x_{1},..,x_{m}])\geq (\varphi -\psi
)([1,\zeta_{0}^{[n-1]},\zeta_{1}^{[n]},..,\zeta_{k-1}^{[n]}]),
\end{eqnarray}
o\`u $\zeta^{[d]}=(\zeta ,..,\zeta)\in \complex^{d}$,
$\zeta_{0}=(x_{1}..x_{n-1})^{1/(n-1)}$, et pour $1\leq h\leq k-1$,
$\zeta_{h}=(x_{hn}..x_{(h+1)n-1})^{1/n}$.
\end{lem}
{\bf Preuve.} La d\'emonstration s'effectue ind\'ependemment sur
chacun des $n$-uplets $(x_{hn},..,x_{(h+1)n-1})\in \reel^{n}$,
avec $x_{0}=1$. En effet, sous les hypoth\`eses du lemme
\ref{lem1}, l'in\'egalit\'e (\ref{eq1}) est la cons\'equence de
l'in\'egalit\'e:
\begin{eqnarray}\label{eq1-0}
(\varphi -\psi )([1,x_{1},..,x_{m}])\geq (\varphi -\psi
)([1,\zeta_{0}^{[n-1]},x_{n},..,x_{m}]),
\end{eqnarray}
et des $k-1$ in\'egalit\'es:
\begin{eqnarray}\label{eq1-h}
(\varphi -\psi )([1,x_{1},..,x_{m}])\geq (\varphi -\psi
)([1,x_{1},..,x_{(h-1)n-1};\zeta_{h}^{[n]};x_{(h+1)n},..,x_{m}]),
\end{eqnarray}
pour $1<h \leq (n-1)$ .

La preuve de ces $k$ in\'egalit\'es \'etants identiques, on se
propose de r\'ediger celle de l'in\'egalit\'e (\ref{eq1-0}). La
d\'emonstration se fait par r\'ecurrence. Supposons que l'on ait
pour tout $(x_{1},..,x_{m})$ avec $0< x_{i}\leq 1$ et pour $1\leq
p < n-1$ :
\begin{eqnarray}\label{eq2}
&&(\varphi -\psi )([1,x_{1},..,x_{m}])\geq  \nonumber \\
&&(\varphi -\psi
)([1,(x_{1}..x_{p})^{1/p},..,(x_{1}..x_{p})^{1/p},
x_{p+1},...,x_{n-1};x_{n},...,x_{m}]),
\end{eqnarray}
l'hypoth\`ese de r\'ecurrence \'etant  \'evidente pour $p=1$.
Supposons qu'elle n'est pas satisfaite au rang $p+1$. Il
existerait alors un point $(x_{1}^{0},...,x_{m}^{0})\in \reel^{m}$
tel que
\begin{eqnarray}\label{eq3}
&&(\varphi -\psi )([1,x_{1}^{0},..,x_{m}^{0}])< (\varphi
-\psi )([1,(x_{1}^{0}..x_{p+1}^{0})^{1/(p+1)},..,\nonumber \\
&&(x_{1}^{0}..x_{p+1}^{0})^{1/(p+1)},
x_{p+2}^{0},...,x_{n-1}^{0};x_{n}^{0},...,x_{m}^{0}]),
\end{eqnarray}
En utilisant la $G_{n,k}$-invariance de $\varphi$, on peut
supposer que $x_{1}^{0}\leq ...\leq x_{p+1}^{0}$. D'autre part, en
appliquant la $G_{n,k}$-invariance de $\varphi$, et l'hypoth\`ese
de r\'ecurrence (\ref{eq2}) en les points
$$([1,x_{1}^{0},..,x_{p}^{0},
x_{p+1}^{0},...,x_{n-1}^{0};x_{n}^{0},...,x_{m}^{0}])$$ et
$$([1,x_{2}^{0},..,x_{p+1}^{0},
x_{1}^{0},x_{p+2}^{0},...,x_{n-1}^{0};x_{n}^{0},...,x_{m}^{0}]),$$
on a :
\begin{eqnarray}\label{eq2'}
&&(\varphi -\psi )([1,x_{1}^{0},..,x_{m}^{0}])\geq  \nonumber
\\ &&(\varphi -\psi )([1,(x_{1}^{0}..x_{p}^{0})^{1/p},..,
(x_{1}^{0}..x_{p}^{0})^{1/p},
x_{p+1}^{0},...,x_{n-1}^{0};x_{n}^{0},...,x_{m}^{0}]),
\end{eqnarray}
et
\begin{eqnarray}\label{eq2''}
&&(\varphi -\psi )([1,x_{1}^{0},..,x_{m}^{0}])=  \nonumber \\
&&(\varphi -\psi )([1,x_{2}^{0},..,x_{p+1}^{0},
x_{1}^{0},x_{p+2}^{0},...,x_{n-1}^{0};x_{n}^{0},...,x_{m}^{0}])\geq\nonumber
\\ &&(\varphi -\psi
)([1,(x_{2}^{0}..x_{p+1}^{0})^{1/p},..,(x_{2}^{0}..x_{p+1}^{0})^{1/p},
x_{1}^{0} ,...,x_{n-1}^{0};x_{n}^{0},...,x_{m}^{0}]).
\end{eqnarray}
Consid\'erons maintenant la courbe $C$ d'\'equation
$t^{p}x_{p+1}=x_{1}^{0}...x_{p+1}^{0}$ du plan r\'eel :
$$([1,t,...,t,x_{p+1},x_{p+2}^{0},..,x_{m}^{0}]).$$ Les points
$$P_{1}=([1,(x_{1}^{0}..x_{p}^{0})^{1/p},..,
(x_{1}^{0}..x_{p}^{0})^{1/p},
x_{p+1}^{0},...,x_{n-1}^{0};x_{n}^{0},...,x_{m}^{0}])$$ et
$$P_{2}=([1,(x_{2}^{0}..x_{p+1}^{0})^{1/p},..,(x_{2}^{0}..x_{p+1}^{0})^{1/p},
x_{1}^{0} ,...,x_{n-1}^{0};x_{n}^{0},...,x_{m}^{0}])$$
appartiennent \`a l'\'evidence \`a cette courbe. Si
$x_{1}^{0}=...=x_{p+1}^{0}$, il n'y a rien \`a montrer. Sinon,
sachant que l'on a choisi $x_{1}^{0}\leq ...\leq x_{p+1}^{0}$, les
points $P_{1}$ et $P_{2}$, se trouvent de part et d'autre de la
diagonale $t=x_{p+1}$. Or la courbe $C$ intersecte cette diagonale
en le point
$$P_{3}=([1,(x_{1}^{0}..x_{p+1}^{0})^{1/(p+1)},..,(x_{1}^{0}..x_{p+1}^{0})^{1/(p+1)},
x_{p+2}^{0},...,x_{n-1}^{0};x_{n}^{0},...,x_{m}^{0}])$$
intervenant dans la relation (\ref{eq3}). D'autre part, en
utilisant les relations (\ref{eq3}), (\ref{eq2'}) et (\ref{eq2''})
on a : $(\varphi - \psi )(P_{3})\geq (\varphi - \psi
)(P_{1})$ et $(\varphi - \psi )(P_{3})\geq (\varphi - \psi
)(P_{2})$, ce qui prouve l'existence d'un maximum local pour la
fonction $(\varphi - \psi )$ sur la courbe $C$. Ceci implique,
au vu de la $G_{n,k}$-invariance de $(\varphi - \psi )$
l'existence d'un maximum local en un point $P$ pour cette fonction
sur la courbe holomorphe (toujours not\'ee $C$) d'\'equation
$\zeta^{p}z_{p+1}=x_{1}^{0}...x_{p+1}^{0}$ du plan complexe
$$([1,\zeta,...,\zeta,z_{p+1},x_{p+2}^{0},..,x_{m}^{0}]).$$ La
Hessienne complexe de terme g\'{e}n\'{e}ral $$\frac
{\partial^{2}}{\partial z
\partial \overline{z}} \{(\varphi
-\psi)(C(z))\} =\frac{\partial^{2}(\varphi
-\psi)}{\partial z_{\lambda}
\partial \overline{z}_{\mu} }(C(z))C^{\lambda}(z)\overline{C}^{\mu}(z)$$
est alors d\'{e}finie n\'{e}gative en $P$.

D'autre part, sachant que $-\frac{\partial^{2}\psi}{\partial
z_{\lambda}
\partial \overline{z}_{\mu}}= g_{\lambda\overline{\mu}},$
$$g_{\lambda\overline{\mu}}+\frac{\partial^{2}\varphi}{\partial
z_{\lambda}
\partial \overline{z}_{\mu}}=\frac{\partial^{2}(\varphi -\psi )}{\partial z_{\lambda}
\partial \overline{z}_{\mu}}$$ sera donc, d'apr\`es ce qui pr\'ec\'ede, d\'efinie n\'egative en $P$,
ce qui est en contradiction avec la $g$-admissibilit\'{e} de
$\varphi$. D'o\`u l'in\'egalit\'e (\ref{eq1-0}) et par
cons\'equent et le lemme \ref{lem1}.


\begin{lem}\label{lem2}
Etant donn\'{e}e une fonction $\varphi\in
C^{\infty}(\proj_{m}\complex)$, $g$-admissible sur
$\proj_{m}\complex$, avec $m=kn-1$, $G_{n,k}$-invariante, on a
pour $0<\zeta_{i}\leq 1$:
\begin{eqnarray}\label{eq4}
(\varphi - \psi
)([1,\zeta_{0}^{[n-1]},\zeta_{1}^{[n]},..,\zeta_{k-1}^{[n]}])\geq
(\varphi - \psi )([1,\zeta_{0}^{[n-1]},\gamma^{[(k-1)n]}]),
\end{eqnarray}
o\`u $\zeta^{[d]}=(\zeta ,..,\zeta)\in \complex^{d}$, et
$\gamma=(\zeta_{1}..\zeta_{k-1})^{1/(k-1)}$.
\end{lem}
{\bf Preuve.} Comme dans le lemme \ref{lem1}, la preuve s'effectue
par r\'ecurrence. Supposons que l'on ait pour tout $0<
\zeta_{i}\leq 1$ et pour $1\leq p < k-1$ :
\begin{eqnarray}\label{eq5}
&&(\varphi -\psi
)([1,\zeta_{0}^{[(n-1)]},\zeta_{1}^{[n]},..,\zeta_{k-1}^{[n]}])\geq
\nonumber \\ &&(\varphi -\psi
)([1,\zeta_{0}^{[(n-1)]},((\zeta_{1}..\zeta_{p})^{1/p})^{[np]},
\zeta_{p+1}^{[n]},...,\zeta_{k-1}^{[n]}]),
\end{eqnarray}
l'hypoth\`ese de r\'ecurrence \'etant  \'evidente pour $p=1$.
Supposons qu'elle n'est pas satisfaite au rang $p+1$. Il
existerait alors un $k$-uplet $(\delta_{0},...,\delta_{k-1})\in
\reel^{k}$ tel que
\begin{eqnarray}\label{eq6}
&&(\varphi -\psi
)([1,\delta_{0}^{[(n-1)]},\delta_{1}^{[n]},..,\delta_{k-1}^{[n]}])
<  \nonumber \\ &&(\varphi -\psi
)([1,\delta_{0}^{[(n-1)]},((\delta_{1}..\delta_{p+1})^{1/(
p+1)})^{[n(p+1)]}, \delta_{p+2}^{[n]},...,\delta_{k-1}^{[n]}]),
\end{eqnarray}
Quitte \`a \'echanger l'emplacement des $n$-uplets
$\delta_{i}^{[n]}$ en utilisant la $G_{n,k}$-invariance de
$\varphi$, on peut supposer que $\delta_{1}\leq ...\leq
\delta_{p+1}$. Puis en appliquant, comme dans le lemme
pr\'ec\'edant, l'hypoth\`ese de r\'ecurrence (\ref{eq5}) en les
points
$$([1,\delta_{0}^{[(n-1)]},\delta_{1}^{[n]},..,\delta_{p}^{[n]},\delta_{p+1}^{[n]},..,\delta_{k-1}^{[n]}])$$
et
$$([1,\delta_{0}^{[(n-1)]},\delta_{2}^{[n]},..,\delta_{p+1}^{[n]},\delta_{1}^{[n]},
\delta_{p+2}^{[n]},..,\delta_{k-1}^{[n]}]),$$ et en utilisant
encore la $G_{n,k}$-invariance de $\varphi$, on a :
\begin{eqnarray}\label{eq5'}
&&(\varphi -\psi
)([1,\delta_{0}^{[(n-1)]},\delta_{1}^{[n]},..,\delta_{k-1}^{[n]}])\geq
\nonumber \\ &&(\varphi -\psi
)([1,\delta_{0}^{[(n-1)]},((\delta_{1}..\delta_{p})^{1/p})^{[np]},
\delta_{p+1}^{[n]},...,\delta_{k-1}^{[n]}]),
\end{eqnarray}
et
\begin{eqnarray}\label{eq5''}
&&(\varphi -\psi
)([1,\delta_{0}^{[(n-1)]},\delta_{1}^{[n]},..,\delta_{k-1}^{[n]}])\geq
\nonumber \\ &&(\varphi -\psi
)([1,\delta_{0}^{[(n-1)]},((\delta_{2}..\delta_{p+1})^{1/p})^{[np]},
\delta_{1}^{[n]}, \delta_{p+2}^{[n]},...,\delta_{k-1}^{[n]}]),
\end{eqnarray}
Consid\'erons maintenant la courbe $\Gamma$ d'\'equation
$\zeta^{p}\zeta_{p+1} =\delta_{1}...\delta_{p+1}$ (on rappelle que
les $\delta$ sont fixes) du plan r\'eel :
$$([1,\delta_{0}^{[(n-1])},\zeta^{[np]},\zeta^{[n]}_{p+1},\delta^{[n]}_{p+2},..,\delta^{[n]}_{k-1}]).$$
Les points
$$Q_{1}=([1,\delta_{0}^{[(n-1)]},((\delta_{1}..\delta_{p})^{1/p})^{[np]},
\delta_{p+1}^{[n]},...,\delta_{k-1}^{[n]}])$$ et
$$Q_{2}=([1,\delta_{0}^{[(n-1)]},((\delta_{2}..\delta_{p+1})^{1/p})^{[np]},
\delta_{1}^{[n]}, \delta_{p+2}^{[n]},...,\delta_{k-1}^{[n]}])$$
appartiennent \`a $\Gamma$. Si $\delta_{1}=...=\delta_{p+1}$, il
n'y a rien \`a montrer. Sinon, par le choix $\delta_{1}\leq
...\leq \delta_{p+1}$ fait plus haut, les points $Q_{1}$ et
$Q_{2}$ se trouvent de part et d'autre de la diagonale $\zeta
=\zeta_{p+1}$ qui intersecte la courbe $\Gamma$ en le point
$$Q_{3}=([1,\delta_{0}^{[(n-1)]},((\delta_{1}..\delta_{p+1})^{1/(
p+1)})^{[n(p+1)]}, \delta_{p+2}^{[n]},...,\delta_{k-1}^{[n]}])$$
de la relation (\ref{eq6}). Les relations (\ref{eq6}),
(\ref{eq5'}) et (\ref{eq5''}) donnent : $(\varphi - \psi
)(Q_{3})\geq (\varphi - \psi )(Q_{1})$ et $(\varphi - \psi
)(Q_{3})\geq (\varphi - \psi )(Q_{2})$, ce qui prouve l'existence
d'un maximum local pour la fonction $(\varphi - \psi )$ sur la
courbe $\Gamma$. $(\varphi - \psi )$ poss\`ede alors un maximum
local en un point $Q$ sur la courbe holomorphe, encore not\'ee
$\Gamma$, d'\'equation $\zeta^{p}\zeta_{p+1}
=\delta_{1}...\delta_{p+1}$ du plan complexe
$$([1,\delta_{0}^{[(n-1])},\zeta^{[np]},\zeta^{[n]}_{p+1},\delta^{[n]}_{p+2},..,\delta^{[n]}_{k-1}]).$$
La Hessienne complexe en $Q$ de terme g\'{e}n\'{e}ral
\begin{eqnarray}
\frac {\partial^{2}}{\partial z
\partial \overline{z}} \{(\varphi
-\psi)(\Gamma (Q))\} =\frac{\partial^{2}(\varphi
-\psi)}{\partial z_{\lambda}
\partial \overline{z}_{\mu} }(\Gamma (Q))\Gamma^{\lambda}(Q)\overline{\Gamma}^{\mu}(Q)\nonumber\\
=
(g_{\lambda\overline{\mu}}+\frac{\partial^{2}\varphi}{\partial
z_{\lambda}
\partial \overline{z}_{\mu}})(\Gamma (Q))\Gamma^{\lambda}(Q)\overline{\Gamma}^{\mu}(Q)\nonumber
\end{eqnarray}
est alors d\'{e}finie n\'{e}gative, ce qui est en contradiction
avec la $g$-admissibilit\'{e} de $\varphi$. D'o\`u
l'in\'egalit\'e (\ref{eq4}) et le lemme \ref{lem2}.


\begin{lem}\label{lem3}
Etant donn\'{e}e une fonction $\varphi\in
C^{\infty}(\proj_{m}\complex)$, $g$-admissible sur
$\proj_{m}\complex$, avec $m=kn-1$, $G_{n,k}$-invariante, on a
pour $0<\zeta , \gamma \leq 1$:
\begin{eqnarray}\label{eq7}
(\varphi - \psi )([1,\zeta^{[n-1]};\gamma^{[(k-1)n]}])\geq
(\varphi - \psi )([1,...,1]),
\end{eqnarray}
o\`u $\zeta^{[d]}=(\zeta ,..,\zeta)\in \complex^{d}$.
\end{lem}
{\bf Preuve.} $(\varphi - \psi )$ \'etant
$G_{n,k}$-invariante, on a
\begin{eqnarray}
& &(\varphi - \psi
)([1,\zeta^{[n-1]};\gamma^{[(k-1)n]}])\nonumber\\ &=&(\varphi -
\psi
)([\gamma^{[n]};1,\zeta^{[n-1]};\gamma^{[(k-2)n]}])\nonumber\\
&=&(\varphi - \psi )([1^{[n]};1/\gamma,(\zeta /\gamma
)^{[n-1]};1^{[(k-2)n]}]).\nonumber
\end{eqnarray}
En appliquant maintenant l'in\'egalit\'e (\ref{eq1-h}) du lemme
\ref{lem1}, au second $n$-uplet de l'\'egalit\'e pr\'ec\'edente,
on a :
\begin{eqnarray}
(\varphi - \psi )([1,\zeta^{[n-1]};\gamma^{[(k-1)n]}]) \geq
(\varphi - \psi
)([1^{[n]};\delta^{[n]};1^{[(k-2)n]}]),\nonumber
\end{eqnarray}
o\`u $\delta = (\zeta^{(n-1)/n})/\gamma$. Enfin, en appliquant le
lemme \ref{lem2}, aux $(k-1)$ derniers $n$-uplets de l'\'egalit\'e
pr\'ec\'edente, on a :
\begin{eqnarray}
(\varphi - \psi )([1,\zeta^{[n-1]};\gamma^{[(k-1)n]}]) \geq
(\varphi - \psi )([1^{[n]};\nu^{[(k-1)n]}]),\nonumber
\end{eqnarray}
o\`u $\nu = \delta^{1/(k-1)}$. Le lemme \ref{lem3} sera donc
\'etabli d\`es que l'on aura prouv\'e, pour tout $\nu >0$,
l'in\'egalit\'e:
\begin{eqnarray}\label{eq7'}
(\varphi - \psi )([1^{[n]};\nu^{[(k-1)n]}])\geq (\varphi -
\psi )([1,...,1]),
\end{eqnarray}
Supposons que l'in\'egalit\'e (\ref{eq7'}) ne soit pas satisfaite
pour un r\'eel $\nu_{0}\neq 1$, on a
\begin{eqnarray}\label{eq7''}
(\varphi - \psi )([1^{[n]};\nu_{0}^{[(k-1)n]}])< (\varphi -
\psi )([1,...,1]),
\end{eqnarray}
En utilisant encore une fois la $G_{n,k}$-invariance de $(\varphi
- \psi )$, et en appliquant le lemme \ref{lem2} aux $(k-1)$
derniers $n$-uplets, on a :
\begin{eqnarray}\label{eq7'''}
&&(\varphi - \psi )([1^{[n]};\nu_{0}^{[(k-1)n]}])\nonumber\\
&=&(\varphi - \psi
)([(1/\nu_{0})^{[n]};1^{[(k-1)n]}])\nonumber\\ &=&(\varphi -
\psi )([1^{[n]};(1/\nu_{0})^{[n]};1^{[(k-2)n]}]),\nonumber\\
&\geq & (\varphi - \psi
)([1^{[n]};\{(1/\nu_{0})^{1/(k-1)}\}^{[(k-1)n]}]).
\end{eqnarray}
Sachant que $\nu_{0}\neq 1$, les r\'eels $\nu_{0}$ et
$(1/\nu_{0})^{1/(k-1)}$ sont situ\'es de part et d'autre du point
$\nu =1$. A l'instar des lemmes pr\'ec\'edents, les in\'egalit\'es
(\ref{eq7''}) et (\ref{eq7'''}) impliquent l'existence d'un
maximum local pour la fonction $(\varphi - \psi )$ sur la
courbe holomorphe d\'ecrite par les points
$[(1^{[n]},z^{[(k-1)n]})]$, ce qui contredit le fait que $\varphi
-\psi$ soit $g$-admissible.


\begin{lem}\label{lem4}
Etant donn\'{e}e une fonction $\varphi\in
C^{\infty}(\proj_{m}\complex)$, $g$-admissible sur
$\proj_{m}\complex$, avec $m=kn-1$, $G_{n,k}$-invariante,
v\'{e}rifiant $\sup\varphi  =0$, sur $\proj_{m}\complex$, on a:\\
$(\varphi -\psi )(1,..,1)\geq 0.$
\end{lem}
{\bf Preuve.} Sachant que $\psi (1,..,1)= -a_{m}\ln (m+1)$, il
suffit de montrer que l'on a $\varphi (1,..,1)\geq -a_{m}\ln
(m+1)$. On raisonne sur la g\'eographie du point $P_{0}$ o\`u
$\sup\varphi$ est atteint. $\varphi$ \'etant $G_{n,k}$-invariante
on peut supposer que $P_{0}=[1,x_{1}^{0},...,x_{m}^{0}]$ du cube
r\'eel $0\leq x_{i}\leq 1$. Remarquons d'abords, que si
$\sup\varphi$ $(=0)$ est atteint en $[1,..,1]$, alors il n'y a
rien \`a montrer.
\begin{itemize}
\item {\underline{Premier cas:} $0< x_{i}^{0}< 1$.}
\end{itemize}
Si $p$ d\'esigne l'indice du plus petit des $x_{i}^{0}$ on
consid\`ere alors la courbe $\Omega$ constitu\'ee des points
 $[1,x_{1},...,x_{m}]$ et v\'erifiant les \'equations :
\begin{eqnarray}
&&x_{1}=x_{p}^{(\ln x_{1}^{0})/(\ln
x_{p}^{0})},..,x_{p-1}=x_{p}^{(\ln x_{p-1}^{0})/(\ln
x_{p}^{0})},\nonumber\\ &&x_{p},x_{p+1}=x_{p}^{(\ln
x_{p+1}^{0})/(\ln x_{p}^{0})},..,x_{m}=x_{p}^{(\ln x_{m}^{0})/(\ln
x_{p}^{0})}.\nonumber
\end{eqnarray}
Cette courbe passe par les points $[1,0,...,0]$, $P_{0}$ et
$[1,..,1]$. Consid\'erons maintenant la fonction $$\tilde\psi
=\ln \frac{\mid z_{0}\mid^{2a_{m}}}{(\mid z_{0}\mid^{2}+...+\mid
z_{m}\mid^{2})^{a_{m}}}.$$ Elle est d\'{e}finie sur
$\complex^{m}-\{z_{0}=0\}$, homog\`{e}ne de degr\'e z\'ero, c'est
donc une fonction sur $\proj_{m}$. $\tilde\psi$ atteint son
maximum (\'egal \`a z\'ero) $[1,0,..,0] \in\proj_{m}\complex$ et
tends vers moins l'infini lorsque la coordonn\'{e}e homog\`{e}ne
$z_{0}$ tend vers z\'{e}ro. Son expression dans la carte
$\{z_{0}=1\}$ est $$\tilde\psi =\ln \frac{1}{(1+\mid
z_{1}\mid^{2}+...+\mid z_{m}\mid^{2})^{a_{m}}}.$$ Supposons que
l'on ait, pour une fonction $\varphi$ donn\'ee, $\varphi
([1,..,1])<-a_{m}\ln (m+1)$ et \'evaluons dans ce cas
$(\varphi-\tilde\psi)$ en les points $[1,0,...,0]$, $P_{0}$ et
$[1,..,1]$ de la courbe $\Omega$. On a :
$$(\varphi-\tilde\psi)([1,0,...,0])=\varphi
([1,0,...,0])\leq0,$$
$$(\varphi-\tilde\psi)(P_{0})=-\tilde\psi (P_{0}) >0,$$ et
enfin $$(\varphi-\tilde\psi)([1,..,1])<-a_{m}\ln (m+1) -
\tilde\psi ([1,..,1])=0.$$ La fonction
$\varphi-\tilde\psi$ poss\`ederait alors un maximum local sur
la courbe holomorphe d\'efinie par $\Omega$, ce qui, encore une
fois mettrait \`a d\'efaut la $g^{m}$-admissibilit\'e de
$\varphi$.
\begin{itemize}
\item {\underline{Second cas:} Un au moins des $x_{i}^{0}$ se
trouve sur une ar\^ete du cube: $0\leq x_{i}\leq 1$, avec les
$x_{i}^{0}$ non tous nuls.}\\
\end{itemize}
Sachant que $$(\varphi-\tilde\psi)(P_{0})=-\tilde\psi
(P_{0}) >0,$$ la continuit\'e de $(\varphi-\tilde\psi)$ permet
d'affirmer l'existence d'un point
$\tilde{P_{0}}=[1,\tilde{x}_{1}^{0},...,\tilde{x}_{m}^{0}]$, avec
$0<\tilde{x}_{i}^{0}<1$, voisin de $P_{0}$ tel que
$$(\varphi-\tilde\psi)(\tilde{P_{0}}) >0,$$ le reste de la
preuve est identique \`a celle du premier cas, en remplacant
$P_{0}$ par $\tilde{P_{0}}$.
\begin{itemize}
\item {\underline{Troisi\`eme cas:} Tous les $x_{i}^{0}$ sont nuls.}
\end{itemize}
On peut encore trouver au voisinage du point $[1,0,..,0]$ un point
$\tilde{P_{0}}=[1,\tilde{x}_{1}^{0},...,\tilde{x}_{m}^{0}]$, avec
$0<\tilde{x}_{i}^{0}<1$ tel que
$$(\varphi-\tilde\psi)(\tilde{P_{0}}) >0.$$ En effet,
supposons qu'il existe un voisinage de $[1,0,..,0]$ ( constitu\'e
de points $[1,x_{1},..,x_{m}]$ dans lequel
$(\varphi-\tilde\psi) \leq 0$. Sachant que
$(\varphi-\tilde\psi)([1,0,..,0]) =0$, la fonction
$\varphi-\tilde\psi$ poss\`ederait dans ce voisinage un
maximum en $[1,0,..,0]$, ce qui contredirait l'admissibilit\'e de
$\varphi$. On peut donc, encore une fois, se ramener au premier
cas, en remplacant $P_{0}$ par
$\tilde{P_{0}}=[1,\tilde{x}_{1}^{0},...,\tilde{x}_{m}^{0}]$.


{\subsection{Preuve du th\'eor\`eme \ref{th1}.}}
Soit $\varphi\in C^{\infty}(\proj_{m}\complex)$, une fonction
$g$-admissible sur $\proj_{m}\complex$, avec $m=kn-1$,
$G_{n,k}$-invariante. Les in\'egalit\'es (\ref{eq1}) et
(\ref{eq4}) des lemmes \ref{lem1} et \ref{lem2} donnent pour,
$0<x_{i}\leq 1$,
\begin{eqnarray}\label{eq8}
(\varphi -\psi )([1,x_{1},..,x_{m}])\geq (\varphi -\psi
)([1,\zeta_{0}^{[n-1]},\gamma^{[(k-1)n]}]),
\end{eqnarray}
o\`u $\zeta^{[d]}=(\zeta ,..,\zeta)\in \complex^{d}$,
$\zeta_{0}=(x_{1}..x_{n-1})^{1/(n-1)}$ et
$$\gamma=\{\prod_{h=1}^{k-1}(x_{hn}..x_{(h+1)n-1})^{1/n}\}^{1/(k-1)}.$$
Puis en utilisant l'in\'egalit\'e (\ref{eq7}) du lemme \ref{lem3}
et le lemme \ref{lem4} , on \'etablit, toujours pour $0<x_{i}\leq
1$,
\begin{eqnarray}\label{eq9}
(\varphi -\psi )([1,x_{1},..,x_{m}])\geq (\varphi -\psi
)([1,...,1])\geq 0.
\end{eqnarray}
$\varphi$ et $\psi$ \'etants $G_{n,k}$-invariantes, la condition
$\varphi\geq\psi$ est v\'erifi\'ee partout d\`es qu'elle l'est en
les points $([1,x_{1},..,x_{m}])$ avec $0\leq x_{i}\leq 1$.
L'in\'egalit\'e (\ref{eq9}) donne le r\'esultat recherch\'e pour
les les points $([1,x_{1},..,x_{m}])$ avec $0< x_{i}\leq 1$. Or
pour les points ayant une composante homog\`ene nulle on a $\psi
=-\infty$, donc \`a l'\'evidence $\varphi\geq\psi$.
L'in\'egalit\'e (\ref{eq9}) et cette derni\`ere remarque
ach\`event donc la preuve du th\'eor\`eme \ref{th1}.

{\subsection{Preuve du th\'eor\`eme \ref{th2}.}}
Soit $\varphi\in C^{\infty}(\proj_{m}\complex )$ une fonction
$g$-admissible et $G$-invariante, v\'{e}rifiant $\sup\varphi = 0$
sur $\proj_{m}\complex$. D'apr\`es le th\'eor\`eme \ref{th1}, on a
: $\varphi \geq \psi$. On a donc pour tout $\alpha$ positif
l'in\'egalit\'e : $$\int_{\proj_{m}\complex} \exp\{-\alpha\varphi
\} dv \leq \int_{\proj_{m}\complex}\exp\{-\alpha\psi \}dv .$$
Calculons cette derni\`ere int\'egrale la carte dense d\'efinie
par $\{z_{0}=1\}$ (cf. introduction). Sachant que l'on a par
hypoth\`ese $a_{m}=m+1$, l'\'el\'ement de volume sera donn\'e,
dans la carte choisie, par : $$dv=(-1)^{m}\frac{dz_{1}\wedge
d\overline{z}_{1}\wedge ...\wedge dz_{m}\wedge d\overline{z}_{m}}
{(1+\mid z_{1}\mid^{2}+...+\mid z_{m}\mid^{2})^{m+1}}.$$ En
utilisant le fait que $\psi$ ne d\'epend que des $\mid z_{p}\mid$,
le changement de variables $x_{p}=\mid z_{p}\mid^{2}$ donne :
\begin{eqnarray}
\int_{\proj_{m}\complex}\exp\{-\alpha\psi \}dv  = Cst
\int_{0}^{+\infty}...\int_{0}^{+\infty}\frac{ (1+
x_{1}+...+x_{m})^{(\alpha -1)(m+1)}dx_{1}...dx_{m}}
{(x_{1}...x_{m})^{\alpha}}\nonumber
\end{eqnarray}
qui converge pour $\alpha < 1$, d'o\`u le r\'esultat, et la preuve
du th\'eor\`eme \ref{th2}.

{\subsection{Preuve du th\'eor\`eme \ref{th3}.}}

Pour la preuve de la premi\`ere partie du th\'eor\`eme \ref{th3},
on se place dans l'ouvert de carte donn\'e par les points
\begin{eqnarray}
&&([\zeta_{0} (z_{0},..,z_{n-1}), \zeta_{1} (z_{n},..,z_{2n-1})
,.., \zeta_{k-1} (z_{(k-1)n},..,z_{kn-1})],\nonumber\\
&&[z_{0},..,z_{n-1}] ,.., [z_{(k-1)n} ,..,z_{kn-1}] )\in
M,\nonumber
\end{eqnarray}
o\`u tous les $\zeta_{i}$ et tous les $z_{i}$ sont non nuls. En
adoptant alors l'\'ecriture non homog\`{e}ne correspondant \`{a}
$z_{0}=1$, l'expression de $\psi_{M}$ (\ref{psi}) dans cette carte
est donn\'{e}e par
\begin{eqnarray}
\psi_{M}=\ln \{ \frac {(1+\mid z_{1}\mid^{2}+..+\mid
z_{kn-1}\mid^{2})^{-k}\mid z_{1}\mid^{2}...\mid
z_{m}\mid^{2}}{(1+\mid z_{1}\mid^{2}+...+\mid
z_{n-1}\mid^{2})^{(n-1)}...(\mid z_{(k-1)n}\mid^{2}+..+\mid
z_{m}\mid^{2})^{(n-1)}}\}.\nonumber
\end{eqnarray}
On a $-\frac{\partial^{2}\psi^{M}}{\partial z_{\lambda}
\partial \overline{z}_{\mu}}= g^{M}_{\lambda\overline{\mu}}.$
En posant, comme pour le th\'{e}or\`{e}me \ref{th1}, $x_{i}=\mid
z_{i}\mid$, les lemmes \ref{lem1}, \ref{lem2} et \ref{lem3}
s'appliquent exactement de la m\^eme mani\`ere en rempla\c{c}ant
$\psi$ par $\psi_{M}$ et les points
$$([1,x_{1},..,x_{n-1};..;x_{(k-1)n},..,x_{m=kn-1}])\in
\mathbb{P}_{m}\mathbb{C}$$ par
$$([1,x_{1},..,x_{m}],[1,x_{1},..,x_{n-1}],..,[x_{(k-1)n},..,x_{m=kn-1}])\in
M.$$ La $G^{M}$-invariance de $\varphi$ et $\psi_{M}$, se
traduisent, dans la carte choisie, de la m\^eme mani\`ere que
celles des fonctions $\varphi$ et $\psi$ d\'efinies dans le
th\'eor\`eme \ref{th1} dans une carte usuelle de
$\mathbb{P}_{m}\mathbb{C}$. On a donc pour $0<x_{i}\leq 1$,
\begin{eqnarray}
&&(\varphi
-\psi_{M})([1,x_{1},..,x_{m}],[x_{1},..,x_{n-1}],..,[x_{(k-1)n},..,x_{m=kn-1}])
\nonumber\\ &&\geq(\varphi -\psi_{M}
)([1,\zeta_{0}^{[n-1]},\gamma^{[(k-1)n]}],[1,\zeta_{0}^{[n-1]}],[\gamma^{[n]}]
,..,[\gamma^{[n]}])\nonumber\\ &&\geq(\varphi -\psi_{M}
)([1^{[m+1]}],[1^{n}],..,[1^{[n]}]) \nonumber
\end{eqnarray}
o\`u $\zeta^{[d]}=(\zeta ,..,\zeta)\in \complex^{d}$,
$\zeta_{0}=(x_{1}..x_{n-1})^{1/(n-1)}$ et
$\gamma=\{\prod_{h=1}^{k-1}(x_{hn}..x_{(h+1)n-1})^{1/n}\}^{1/(k-1)}.$
Pour le lemme \ref{lem4}, on remplace la fonction $\tilde{\psi}$
par $\tilde{\psi}_{M}$ d\'efinie comme suit: On consid\`ere
d'abord la fonction $\breve{\psi}$ d\'efinie sur
$\complex^{m+1}\times(\complex^{n})^{k}
\backslash\bigcup_{i,p}{\{z_{p}^{(i)}=0\}}$, o\`{u} $m=kn-1$, par
\begin{eqnarray}
\breve{\psi} =\ln [\frac{\mid z_{0}^{(0)}\mid^{2k}}{(\mid
z_{0}^{(0)}\mid^{2}+...+\mid z_{m}^{(0)}\mid^{2})^{k}}\times
\frac{\mid z_{0}^{(1)}\mid^{2(n-1)}}{(\mid
z_{0}^{(1)}\mid^{2}+...+\mid
z_{n-1}^{(1)}\mid^{2})^{(n-1)}}\times...\nonumber\\ \times
\frac{\mid z_{0}^{(k)}\mid^{2(n-1)}}{(\mid
z_{0}^{(k)}\mid^{2}+...+\mid
z_{n-1}^{(k)}\mid^{2})^{(n-1)}}].\nonumber
\end{eqnarray}
$(z_{p}^{(0)})_{p}$ \'etant le syst\`eme de coordonn\'ees usuel
sur $\complex^{m+1}$, $(z_{p}^{(1)})_{p}$ celui sur le premier
$\complex^{n}$ du produit,..., et $(z_{p}^{(k)})_{p}$ celui sur le
$k$-i\`{e}me $\complex^{n}$, toujours avec $m=kn-1$).
$\breve{\psi}$ est homog\`{e}ne de degr\'{e} z\'{e}ro sur chacun
des $\complex^{p}$, elle d\'{e}finit par cons\'equent une fonction
(toujours not\'ee $\breve{\psi}$) sur
$\proj_{m}\complex\times\proj_{n-1}\complex\times...\times\proj_{n-1}\complex$
\'egale \`a moins l'infini lorsque la premi\`ere composante
homog\`ene de des $\mathbb{P}_{d}\mathbb{C}$ s'annule. Sa
restriction \`a $M$ sera la fonction $\tilde{\psi}_{M}$. En se
pla\c{c}ant en un point
\begin{eqnarray}
&&([\zeta_{0} (z_{0},..,z_{n-1}), \zeta_{1} (z_{n},..,z_{2n-1})
,.., \zeta_{k-1} (z_{(k-1)n},..,z_{kn-1})],\nonumber\\
&&[z_{0},..,z_{n-1}] ,.., [z_{(k-1)n} ,..,z_{kn-1}] )\in
M,\nonumber
\end{eqnarray}
o\`u tous les $\zeta_{i}$ et tous les $z_{i}$ sont non nuls, on a:
\begin{eqnarray}
\tilde{\psi}_{M}=\ln \{[\mid\zeta_{0}\mid^{2}(\mid
z_{0}\mid^{2}+..+\mid
z_{n-1}\mid^{2})+..+\mid\zeta_{k-1}\mid^{2}(\mid
z_{(k-1)n}\mid^{2}+..+\mid z_{kn-1}\mid^{2})]^{-k}\nonumber\\
\times\frac {\mid \zeta_{0}z_{0}\mid^{2k} \mid
\zeta_{0}z_{0}\mid^{2(n-1)} \mid
\zeta_{1}z_{n}\mid^{2(n-1)}...\mid
\zeta_{k-1}z_{(k-1)n}\mid^{2(n-1)}}{(\mid
\zeta_{0}z_{0}\mid^{2}+...+\mid
\zeta_{0}z_{n-1}\mid^{2})^{(n-1)}...(\mid
\zeta_{k-1}z_{(k-1)n}\mid^{2}+..+\mid
\zeta_{k-1}z_{m}\mid^{2})^{(n-1)}}\},\nonumber
\end{eqnarray}
o\`u $[\zeta_{0},..,\zeta_{k-1}]\in \mathbb{P}_{k}\mathbb{C}$. En
posant pour $i\in\{hn,..,(h+1)n-1\}$, $z'_{i}=\zeta_{h}z_{i}$ on a
\begin{eqnarray}
\tilde{\psi}_{M}=\ln \{\frac {[\mid z'_{0}\mid^{2}+..+\mid
z'_{kn-1}\mid^{2})]^{-k}\mid z'_{0}\mid^{2k} \mid
z'_{0}\mid^{2(n-1)} \mid z'_{n}\mid^{2(n-1)}...\mid
z'_{(k-1)n}\mid^{2(n-1)}}{(\mid z'_{0}\mid^{2}+...+\mid
z'_{n-1}\mid^{2})^{(n-1)}...(\mid z'_{(k-1)n}\mid^{2}+..+\mid
z'_{m}\mid^{2})^{(n-1)}}\}.\nonumber
\end{eqnarray}
Comme dans le lemme \ref{lem4}, on montre, en utilisant des
courbes joignant le point
$([1,0^{[m]}],[1,0^{[n-1]}],..,[1,0^{[n-1]}])$ au point
$([1^{[m+1]}],[1^{[n]}],..,[1^{[n]}])$ et passant par un point en
lequel $\varphi -\tilde{\psi}_{M}>0$, que l'on a $$(\varphi
-\tilde{\psi}_{M})([1^{[m+1]}],[1^{[n]}],..,[1^{[n]}])=(\varphi
-\psi_{M})([1^{[m+1]}],[1^{[n]}],..,[1^{[n]}])\geq 0,$$ d'o\`u le
r\'esultat: $\varphi \geq \tilde{\psi}_{M}$.

La deuxi\`eme partie du th\'eor\`eme \ref{th3} n\'{e}cessite
l'expression du d\'eterminant de la m\'etrique $g^{M}$. D'apr\`es
\cite{BC0} elle est donn\'e par :
\begin{eqnarray}
k^{k-1}\frac{\prod_{h=1}^{k}[(n-1)(1+\mid z_{1}\mid^{2}+...+\mid
z_{kn-1}\mid^{2})+ k(\mid z_{(h-1)n}\mid^{2}+...+\mid
z_{hn-1}\mid^{2})]^{n-1}} {(1+\mid z_{1}\mid^{2}+...+\mid
z_{kn-1}\mid^{2})^{m+1}\prod_{h=1}^{k} (\mid
z_{(h-1)n}\mid^{2}+...+\mid z_{hn-1}\mid^{2})^{n-1}},\nonumber
\end{eqnarray}
o\`u $z_{0}=1$. Cette expression et celle de $\psi_{M}$ permettent
alors de montrer, \`a l'instar du th\'eor\`eme \ref{th1}, que pour
$0<\beta <1$, $\int_{M}\exp\{-\beta{\psi}_{M} \}dv_{M}<+\infty ,$
o\`u $dv_{M}$ est l'\'el\'ement de volume relatif \`a la
m\'etrique $g^{M}$ sur $M$.\\ L'in\'egalit\'e $\varphi \geq \psi$
pour toute fonction $\varphi\in C^{\infty}(M)$,
$G^{M}$-invariante, $g^{M}$-admissible, ayant un $\sup$ \'egal \`a
z\'ero sur $M$, permet alors d'\'ecrire, pour tout $0<\beta <1 $:

$$\int_{M} \exp\{-\beta\varphi \} dv_{M} \leq
\int_{M}\exp\{-\beta{\psi}_{M} \}dv_{M}<+\infty ,$$

Le th\'eor\`eme de Tian permet alors de statuer sur l'existence
d'une m\'{e}trique d'Einstein-K\"{a}hler dans la m\^eme classe de
K\"{a}hler que $g_{M}$, ce qui termine la preuve du th\'eor\`eme
\ref{th3}.



\begin{center}
--------------------------

Universit\'e Pierre et Marie Curie, Paris, France.

e-mail: benabdes@math.jussieu.fr
\end{center}
\end{document}